# VALUES OF INDEFINITE QUADRATIC FORMS
# AT INTEGRAL POINTS AND FLOWS ON SPACES OF LATTICES

## ARMAND BOREL

This mostly expository paper centers on recently proved conjectures in two areas:

A) A conjecture of A. Oppenheim on the values of real indefinite quadratic forms at integral points.

B) Conjectures of Dani, Raghunathan, and Margulis on closures of orbits in spaces of lattices such as $\mathbf{SL}_n(\mathbf{R})/\mathbf{SL}_n(\mathbf{Z})$.

At first sight, A) belongs to analytic number theory and B) belongs to ergodic and Lie theory, and they seem to be quite unrelated. They are discussed together here because of a very interesting connection between the two pointed out by M. S. Raghunathan, namely, a special case of B) yields a proof of A).

The first main goal of this talk is to describe the Oppenheim conjecture and various refinements and to derive them from one statement about closures of orbits in the space of unimodular lattices in $\mathbf{R}^3$ (see Proposition 2 in §2.3). In §§3 and 4 we put this statement in context and describe more general conjectures and results on orbit closures and invariant probability measures on quotients of Lie groups by discrete subgroups. §5 gives some brief comments on the proofs and further developments. §§6, 7, and 8 are devoted to the so-called $S$-arithmetic setting, where we consider products of real and $p$-adic groups. §6 is concerned with a generalized Oppenheim conjecture; §7 with a generalization of the orbit closure theorem proved by M. Ratner [R8]; and §8 with applications to quadratic forms. Since the subject matter of that last section has not been so far discussed elsewhere, we take this opportunity to present proofs, obtained jointly with G. Prasad. Finally, §9 gives the proof of a lemma on symmetric simple Lie algebras, a special case of which is used in §8.

I am glad to thank M. Ratner and G. Prasad for a number of remarks on, and corrections to, a preliminary version of this paper, thanks to which many typos and inaccuracies have been eliminated.

## I. VALUES OF INDEFINITE QUADRATIC FORMS

## 1. The Oppenheim conjecture.

1.1. In the sequel $F$ denotes a *non-degenerate* quadratic form on $\mathbf{R}^n$ which is *indefinite*, i.e. $F(x) = 0$ for some $x \in \mathbf{R}^n - \{0\}$ or equivalently $F(\mathbf{R}^n) = \mathbf{R}$. It may be written

$$F(x) = \sum_{1 \leqq i,\, j \leqq n} f_{ij} x_i x_j \qquad (f_{ij} = f_{ji} \in \mathbf{R},\ \det(f_{ij}) \neq 0).$$

Unless otherwise stated we assume $n \geqq 3$. We are concerned with $F(\mathbf{Z}^n)$.

**Definition.** $F$ is said to be *rational* if $F(x)/F(y) \in \mathbf{Q}$ whenever $x$, $y \in \mathbf{Q}^n$ and $F(y) \neq 0$, and *irrational* otherwise.

$F$ is rational if and only if there exists $c \in \mathbf{R}^*$ such that $F = c.F_o$, where $F_o$ has rational coefficients (with respect to a basis of $\mathbf{Q}^n$). We may then also arrange $F_o$ to have integral coefficients. Therefore

$$F(\mathbf{Z}^n) = c.F_o(\mathbf{Z}^n) \subset c.\mathbf{Z}$$

is discrete. The Oppenheim conjecture states that, conversely, if $F$ is irrational, then $F(\mathbf{Z}^n)$ is not discrete around the origin. More precisely, consider the two conditions:

(i) *$F$ is irrational.*

(ii) *Given $\epsilon > 0$, there exists $x \in \mathbf{Z}^n$ such that $0 < |F(x)| < \epsilon$.*

We just saw that (ii) $\Rightarrow$ (i). The *Oppenheim conjecture is that* (i) $\Rightarrow$ (ii).

1.2. Historically, this is a bit of an oversimplification. In 1929, A. Oppenheim stated that the following is very likely to be true: if $F$ is irrational and $n \geqq 5$, then $|F(x)|$ takes arbitrary small values on $\mathbf{Z}^n$ [O1, O2]. Formally, this may be written

(ii)′ Given $\epsilon > 0$, there exists $x \in \mathbf{Z}^n$ such that $|F(x)| < \epsilon$,

and is automatically satisfied if $F$ "represents zero rationally"; i.e. if there exists $x \in \mathbf{Q}^n - \{0\}$, hence also $x \in \mathbf{Z}^n - \{0\}$, such that $F(x) = 0$. But Oppenheim had clearly (ii) in mind, and he made it explicit in [O3], still for $n \geqq 5$, though; but [O4, O5] show that he was wondering whether it might be true for $n \geqq 3$ already (it was well known to be false for $n = 2$; an example is given in (1.4)). Then, later, the conjecture (i) $\Rightarrow$ (ii)′ became erroneously known as "Davenport's conjecture", though Davenport referred to Oppenheim. The implication (i) $\Rightarrow$ (ii)′ is obviously equivalent to the Oppenheim conjecture (i) $\Rightarrow$ (ii) for forms not representing zero rationally. The bound $n \geqq 5$ had been suggested to A. Oppenheim by a theorem of A. Meyer according to which a rational indefinite quadratic form in $n \geqq 5$ variables always represents zero rationally (see e.g. [S, IV, §3]). He felt that, in the irrational case, it should take values close to zero on $\mathbf{Z}^n$.

1.3. The condition (ii) leaves open the possibility that $F(\mathbf{Z}^n)$ accumulates to zero only on one side, but Oppenheim showed in [O3] that this cannot happen for $n \geqq 3$, our standing assumption (but that it can for $n = 2$). It then follows by a very elementary argument that (ii) implies:

(iii) $F(\mathbf{Z}^n)$ is dense in $\mathbf{R}$,

so that the conjectural dichotomy was in fact

(A) $\quad$ $F$ *rational* $\Leftrightarrow$ $F(\mathbf{Z}^n)$ *discrete*;

$\quad\quad\quad$ $F$ *irrational* $\Leftrightarrow F(\mathbf{Z}^n)$ *dense.*

1.4. To conclude this section, we give a simple counterexample for $n = 2$, borrowed from [G].

Let $F(x, y) = y^2 - \theta^2.x^2$, where $\theta$ is quadratic, irrational $> 0$, and $\theta^2$ is irrational. The form $F$ is irrational. As is well known, there exists $c > 0$ such that

(1) $\quad\quad\quad\quad |\theta - y/x| \geqq c.x^{-2} \quad\quad (x, y \in \mathbf{Z}, x \neq 0).$

For $x \neq 0$, we can write

$$(2) \qquad F(x,\ y) = x^2(y/x + \theta)(y/x - \theta).$$

We have to prove that $|F(x,\ y)|$ has a strictly positive lower bound for $x,\ y \in \mathbf{Z}$ not both zero. This is clear if one of them is equal to zero. So let $x,\ y \neq 0$. We may assume them to be $> 0$. Then $|\theta + y/x| \geqq \theta$. Together with (1) and (2), this yields $|F(x,\ y)| \geqq c.\theta$.

## 2. Results.

2.1. The first partial results on the Oppenheim conjecture were obtained in the framework of analytic number theory. It was shown to be true for diagonal forms in $n \geqq 9$ variables [C], in $n \geqq 5$ variables [DH] and for general forms in $n \geqq 21$ variables [DR]. Oppenheim himself proved it when $F$ represents zero rationally, for $n \geqq 5$ in [O4] and for $n = 4$ in [O5]. In both papers he stated his belief it should be true for $n = 3$, though it was obviously false for $n = 2$.

It is easy to see that if $F$ is irrational, then there exists a three-dimensional subspace $V \subset \mathbf{Q}^n$ such that the restriction of $F$ to $V \otimes_{\mathbf{Q}} \mathbf{R}$ is non-degenerate, indefinite and irrational. Therefore it suffices to prove the Oppenheim conjecture for $n = 3$.

2.2. Around 1980, M. S. Raghunathan made a conjecture on closures of orbits in spaces of lattices (see 3.2). He noticed further that a very special case, namely Proposition 1 below, would readily imply the "Davenport conjecture". G. A. Margulis, following this strategy, then proved Proposition 1 and deduced from it that (i) $\Rightarrow$ (ii)$'$ [M1, M2]. When informed (by the author, October 1987) of the fact that the Oppenheim conjecture was a slightly stronger one, he quickly completed his argument and established:

**Theorem 1** (Margulis [M3]). *The Oppenheim conjecture is true.*

2.3. The statement on closures of orbits proved and used by Margulis to show that (i) $\Rightarrow$ (ii)$'$ is:

**Proposition 1.** *Let $n = 3$. Then any relatively compact orbit of $SO(F)$ on $\Omega_3 = \mathbf{SL}_3(\mathbf{R})/\mathbf{SL}_3(\mathbf{Z})$ is compact.*

It is of course equivalent to the same assertion for the topological identity component $SO(F)^o$ of $SO(F)$ (which has index two).

A slight extension of it (see Theorem 1$'$ in [M3]) allowed him to establish the full Oppenheim conjecture. Rather than explaining how, I shall sketch a derivation of this conjecture from the following assertion, stronger than Proposition 1 (but still a very special case of the topological conjecture; see §3), proved shortly afterwards for this purpose by S. G. Dani and G. A. Margulis [DM1]:

**Proposition 2.** *We keep the previous notation and assumptions. Then any orbit of $SO(F)^o$ on $\Omega_3$ either is closed and carries an $SO(F)^o$-invariant probability measure or is dense.*

(In this paper, all measures are Borel measures.)

They deduced directly from it that (i) implies not only (iii), but also that the set of values of $F$ on the *primitive vectors* in $\mathbf{Z}^n$ is dense [DM1]. Recall that $x \in \mathbf{Z}^n$, $x \neq 0$, is *primitive* if it is not properly divisible in $\mathbf{Z}^n$. The set of primitive vectors is $\mathbf{SL}_n(\mathbf{Z}).e_1$, where $e_1$ is the first canonical basis vector. This is a sharpening of (iii) which, as far as I know, had never been considered before.

Note that, since $\Omega_3$ is not compact, Proposition 2 obviously implies Proposition 1.

2.4. We now sketch the proof of the Oppenheim conjecture, or rather of its strengthening just mentioned, using Proposition 2.

We let $G = \mathbf{SL}_3(\mathbf{R})$, $H = SO(F)^o$ and $\Gamma = \mathbf{SL}_3(\mathbf{Z})$.

Let $o$ be the origin in $\Omega_3$, i.e. the coset $\Gamma$. Then by Proposition 2, $H.o$ is either closed with finite invariant measure or dense. Assume first it is dense. Then its inverse image $H.\Gamma$ in $G$ is also dense. Fix $c \in \mathbf{R}$. There exists $x \in \mathbf{R}^n$, $x \neq 0$ such that $F(x) = c$. Let $e_1, \ldots, e_n$ be the canonical basis of $\mathbf{R}^n$. There exists $g \in \mathbf{SL}_3(\mathbf{R})$ such that $g.e_1 = x$. Since $H.\Gamma$ is dense in $G$, we can find sequences $\gamma_j \in \Gamma$ and $h_j \in H$ $(j = 1, \ldots)$ such that $h_j.\gamma_j \to g$. Then we have

$$c = F(g.e_1) = \lim_{j \to \infty} F(h_j.\gamma_j.e_1) = \lim_{j \to \infty} F(\gamma_j.e_1);$$

therefore $c$ is in the closure of $F(\mathbf{SL}_3(\mathbf{Z}).e_1)$.

Assume now that $H.o$ is closed and supports an $H$-invariant probability measure. Let $\Gamma_o = H \cap \Gamma$. Then $H/\Gamma_o$ is homeomorphic to $H.o$ and has therefore finite invariant volume. By a general result, this implies that $\Gamma_o$ is "Zariski-dense" in $H$, i.e. is not contained in any algebraic subgroup. However, in this case, it can be checked more directly. In fact, the only important property of real algebraic groups relevant here is that they have only finitely many connected components, in the usual topology. It suffices therefore to show that $\Gamma_o$ is not contained in any closed subgroup of $H$ having finitely many connected components. Let $M$ be one. Then $H/M$ also carries an invariant probability measure (define the measure of an open set $U$ in $H/M$ as equal to that of its inverse image in $H/\Gamma_o$ under the natural projection $H/\Gamma_o \to H/M$). By a standard fact (see e.g. [Bu], VII, §2, no. 6), this implies that $M$, hence also $M^o$, is unimodular. Also, since $H/M^o$ is a finite covering of $H/M$, it carries an invariant probability measure, too, and the same is true if we replace $H$ by its twofold covering $\mathbf{SL}_2(\mathbf{R})$ and $M^o$ by the identity component of its inverse image there. Being unimodular, $M^o$ is not a maximal connected solvable subgroup of $H$. It is then either conjugate to the subgroup $A$ of diagonal matrices with positive entries or to the subgroup $N$ of upper triangular unipotent matrices. Recall the Iwasawa decomposition $H = K.A.N$, where $K = \mathbf{SO}(2)$. If $M^o = A$, then there is an $N$-equivariant diffeomorphism of $H/M^o$ onto $N \times K$. If $x \in N$, $x \neq 1$, then we can find a neighborhood of 1 in $N \times K$, the translates of which by the powers of $x$ are disjoint, hence $H/M^o$ has infinite invariant measure. If $M^o = N$, then there is similarly an $A$-equivariant diffeomorphism of $H/M^o$ onto $A \times K$, and we see in the same way that an invariant measure has infinite volume. Thus, the existence of $M$ leads to a contradiction, which shows that $\Gamma \cap H$ is Zariski-dense in $H$. (In both cases, it would suffice in fact to note that $H/M^o$ is not compact, in view of a theorem of G. D. Mostow [M, Theorem 7.1] which states that if the quotient of a connected Lie group by a closed subgroup with finitely many connected components carries an invariant measure, it is compact. However this result is closely related to the Zariski-density theorem I was trying to prove directly.)

Consider the natural action of $H$ on the space of $3 \times 3$ real symmetric matrices. The only invariants of $H$ are the multiples $cF$ of $F$ $(c \in \mathbf{R})$. Since $\Gamma_o$ is Zariski-dense in $H$, these are also the only invariants of $\Gamma_o$. But invariance under $\Gamma_o$ translates into a system of linear equations for the coefficients of $F$, with integral coefficients. It has therefore a rational solution; i.e. $cF$ has rational coefficients for some $c \neq 0$, hence $F$ is rational by definition.

2.5. The results of M. Ratner recalled in §3 imply that Proposition 2 is valid for all $n \geqq 3$, (see 3.6). As we saw, the case $n = 3$ suffices for the results on quadratic forms discussed so far, but Proposition 2 for arbitrary $n$ also yields a stronger approximation theorem for quadratic forms. To state it, let us say that a subset $(x_1, \dots, x_m)$ of $\mathbf{Z}^n$ $(m \leqq n)$ is primitive if it is part of a basis of $\mathbf{Z}^n$. If $m < n$, this condition is equivalent to the existence of $g \in \mathbf{SL}_n(\mathbf{Z})$ such that $g.e_i = x_i$ $(i = 1, \dots, m)$. Then we have [BP2, 7.9].

**Theorem 2.** *Let $c_i \in \mathbf{R}$ $(i = 1, \dots, n-1)$. Assume $F$ to be irrational. Then there exists a sequence $(x_{j1}, \dots, x_{jn-1})$ $(j = 1, 2, \dots)$ of primitive subsets of $\mathbf{Z}^n$ such that*

$$(1) \qquad \lim_{j \to \infty} F(x_{ji}) = c_i \qquad (i = 1, \dots, n-1).$$

Proposition 2 already implies it for two values $c_1$, $c_2$, as was shown in [DM1]. The proof is an easy extension of the first argument in 2.4. In fact [BP, 7.9] is concerned more generally with the $S$-arithmetic case (see §6) but without finite places. The proof in the general $S$-arithmetic case will be given in 8.4.

2.6. Propositions 1 and 2 are very special cases of the general results of M. Ratner outlined in the next section. However, Proposition 1 had earlier been given a comparatively elementary proof in [M1, M3], and [DM3] provides also an elementary proof of a theorem weaker than Proposition 2 but stronger than those of [M1, M3] and already sufficient to show that if $F$ is irrational, its set of values at primitive vectors is dense in $\mathbf{R}$. The main theorems on flows, in particular the property of "uniform distribution" (see 3.9), have led to some quantitative refinements of the Oppenheim conjecture. We state here the simplest one

**Proposition 3.** *Given $a$, $b > 0$ with $a < b$, there exists constants $r_o$, $c > 0$ such that*

$$Card\{x \in \mathbf{Z}^n \cap B_r \big| a \leqq |F(x)| \leqq b\} \geqq c.r^{n-2} \qquad (r \geqq r_o).$$

Here, $B_r$ denotes the euclidean ball in $\mathbf{R}^n$ of radius $r$ with center the origin. This was proved first, independently, by S. G. Dani and S. Mozes on the one hand and M. Ratner on the other (unpublished). A more general statement, valid for certain compact sets of quadratic forms, is proved in [DM4]; see Corollary on page 95.

2.7. The truth of the Oppenheim conjecture also yields a characterization of arbitrary non-degenerate rational quadratic forms in $n \geqq 2$ variables, as was noted by G. Prasad and me [B]:

**Proposition 4.** *Let $n \geqq 2$ and $E$ be a non-degenerate quadratic form on $\mathbf{R}^n$. Then $E$ is irrational if and only, given $\epsilon > 0$, there exists $x, y \in \mathbf{Z}^n$ such that*

$$(1) \qquad 0 < |E(x) - E(y)| < \epsilon.$$

This follows by applying Theorem 1 to $E \oplus -E$ on $\mathbf{R}^{2n}$. In fact, in view of 2.4, we may also find primitive vectors $x$, $y$ satisfying (1). Let now $E$ be *positive* non-degenerate. Then $E(\mathbf{Z}^n)$ is discrete in $\mathbf{R}$, but one can still ask questions about the difference between two consecutive values. (1) shows that the lower bound of the non-zero differences $|F(x) - F(y)|$ is zero. This is a small step towards a conjecture made by D. J. Lewis [L], namely: given $\varepsilon > 0$, there exists $R(\varepsilon) > 0$ such that if the norm of $x$ is $\geqq R(\varepsilon)$, then there exists $y \in \mathbf{Z}^n$ so that (1) is satisfied. In other words, if we go far enough, these successive differences are uniformly bounded and tend to zero.

## II. Flows on spaces of lattices

*In this part, $G$ is a connected Lie group and $\Gamma$ is a discrete subgroup which, unless otherwise stated, has finite invariant covolume; i.e. $\Omega = G/\Gamma$ carries a $G$-invariant probability measure, and $H$ is a closed subgroup of $G$.*

*If $L$ is a Lie group, $L^o$ denotes the connected component of the identity in $L$. If $L$ is a group operating on a space $X$, and $x \in X$, then $L_x$ is the isotropy group of $x$, i.e. the subgroup of all elements of $L$ leaving $x$ fixed.*

We recall that the assumption on $\Omega$ forces $G$ to be unimodular. The orbits of $H$ on $\Omega$ define a foliation. We are concerned with the closures of the leaves and the supports of $H$-invariant ergodic probability measures for certain classes of subgroups $H$.

The most important case here is the one where

(MC)                $G = \mathbf{SL}_n(\mathbf{R})$, $\Gamma = \mathbf{SL}_n(\mathbf{Z})$ and $H = SO(F)^o$,

to be referred to as our main special case (MC). Then $\Omega = G/\Gamma$ may be identified with the space of unimodular lattices in $\mathbf{R}^n$.

## 3. The topological conjecture.

3.1. The prototype here is the *horocycle flow* on a Riemann surface of finite area. Let then $G = \mathbf{PSL}_2(\mathbf{R})$ and $X$ be the upper half-plane. Assume, to avoid ramification, that $\Gamma$ is torsion free. Then $X/\Gamma$ is a compact Riemann surface or the complement of finitely many points in one, and $G/\Gamma$ may be identified to the unit tangent bundle of $X/\Gamma$. Take for $H$ the group of matrices

$$\begin{pmatrix} 1 & c \\ 0 & 1 \end{pmatrix} \qquad (c \in \mathbf{R}).$$

Then the orbits of $H$ in $G/\Gamma$ are the orbits of the horocycle flow. By results of Hedlund [H] any orbit of $H$ in $G/\Gamma$ is either compact or dense, and the former does not occur if $G/\Gamma$ is compact.

In [D2], S. G. Dani states a far-reaching conjectural generalization of the previous theorem, proposed by M. S. Raghunathan (Conjecture II, p. 358).

3.2. **Conjecture** (M.S. Raghunathan). *Assume $G$ to be reductive (as a Lie group, i.e. $\mathfrak{g}$ reductive). Let $H$ be an Ad-unipotent (see below) one-parameter subgroup of $G$ and $x \in \Omega$. Then there exists a connected closed subgroup $L$ of $G$ containing $H$ such that $\overline{H.x} = L.x$.*

In [M1, M3], Margulis extended the conjecture to the case where $G$ is a connected Lie group and $H$ is a subgroup generated by elements which are Ad-unipotent in $G$.

An element $g \in G$ is Ad-unipotent if its image Ad $g$ in the adjoint representation of $G$ in its Lie algebra $\mathfrak{g}$ is unipotent (all eigenvalues equal to one). Similarly, $y \in \mathfrak{g}$ is ad-nilpotent if ad $y$ is a nilpotent endomorphism of $\mathfrak{g}$. If $G$ is semi-simple, linear, then $g \in G$ is Ad-unipotent if it is either unipotent as a matrix or central in $G$, and $y \in \mathfrak{g}$ is ad-nilpotent if and only if it is a nilpotent matrix. If $y \in \mathfrak{g}$, then the one-parameter subgroup exp $\mathbf{R}.y$ is Ad-unipotent (i.e. consists of Ad-unipotent elements) if and only if $y$ is ad-nilpotent.

Various special cases of this conjecture were established, notably for horospherical subgroups of reductive groups [D3], until M. Ratner proved it in full generality or, rather, obtained a stronger conclusion under more general assumptions:

3.3. **Theorem 3** (M. Ratner [R4]). *Let $H^o$ be the identity component of $H$. Assume that $H/H^o$ is finitely generated, $H^o$ is generated by Ad-unipotent elements and each coset $h.H^o$ $(h \in H)$ contains an Ad-unipotent element. Let $x \in \Omega$. Then there exists a closed subgroup $L$ of $G$ containing $H$ such that $\overline{H.x} = L.x$ and $L.x$ supports a $L$-invariant probability measure ergodic for the action of $H$.*

*Remark.* The subgroup $L$ is not necessarily unique. For instance, a bigger closed subgroup $L'$ of the same dimension, such that $L'/(L' \cap G_x) = L/(L \cap G_x)$, would also do. However, the Lie algebra $\mathfrak{l}$ of $L$ is uniquely determined: the differential $\mu_1$ at 1 of the map $g \mapsto g.x$ is an isomorphism of $\mathfrak{g}$ onto the tangent space to $\Omega$ at $x$, and $\mathfrak{l}$ is the subspace of $\underline{\mathfrak{g}}$ mapped by $\mu_1$ onto the tangent space to $\overline{H.x}$ at $x$ (which is well defined, since $\overline{H.x}$ is a submanifold by the theorem). Therefore $L^o.H$ is the smallest possible choice for $L$ and is unique. It will be denoted $L(x, H)$. If $H$ is connected, then $L(x, H)$ is the only connected, subgroup satisfying the conclusion of the theorem.

This normalization is introduced in [R3, p. 546].

As an obvious consequence of the above, we have the

3.4. **Corollary.** *Assume that $H$ is connected and maximal among proper closed connected subgroups of $G$. Then any orbit of $H$ in $\Omega$ is either closed and supports an $H$-invariant probability measure or dense. In particular, if $\Omega$ is not compact, a relatively compact orbit is closed and supports an invariant probability measure.*

3.5. The main point in [R4] is to show 3.3 for $H$ connected, one-dimensional, in which case Ratner proves a stronger result (see 3.9), and the following complement:

(∗) *For a given $x \in \Omega$, the set of subgroups $L(x, H)$ occurring in 3.3, when $H$ runs through all the connected closed Ad-unipotent one-dimensional subgroups, is countable.*

3.6. In this section we consider our main special case (MC). Since $F$ is indefinite and $n \geqq 3$, the group $H$ is generated by connected one-dimensional unipotent subgroups. Moreover $SO(F)$ is the fixed point set of an involutive automorphism of $G$; therefore, as is well known and easy to prove, $H$ is maximal among proper closed connected subgroups, so that 3.4 holds and yields in particular Proposition 2 in any dimension $n \geqq 3$.

We now give some indications of how to prove (∗) and reduce the proof of 3.3 for connected $H$'s to one-dimensional ones, in the case under consideration.

We note first (this is completely general) that if 3.3 and 3.5 are true for $H$ and $x$, then they also hold for $g.H.g^{-1}$ and $g.x(g \in G)$. We may therefore assume $x$ to be the origin $o$ of $\Omega$.

(a) If $L.o$ is closed and carries an invariant probability measure, then $L \cap \Gamma$ is of finite covolume in $L$. This implies, rather easily, that $L$ is defined over $\mathbf{Q}$ (see Proposition 1.1 in [BP2]). However $\mathbf{SL}_n(\mathbf{R})$ contains only countably many real algebraic subgroups defined over $\mathbf{Q}$, whence $(*)$.

(b) We now assume 3.3 to be true for all one-dimensional connected Ad-unipotent subgroups of $H$. We want to prove 3.3 for $H$.

Let $\mathcal{N}$ be the variety of nilpotent elements in the Lie algebra $\mathfrak{h}$ of $H$. For $y \in \mathcal{N}$, we let $U_y$ denote the unipotent (or, equivalently, Ad-unipotent) one-dimensional subgroup $\exp \mathbf{R}.y$.

$\mathcal{N}$ is an algebraic variety, invariant under Ad $H$. It is a finite union of irreducible subvarieties $V_j$ $(j \in J)$. Since $\mathcal{N}$ is invariant under Ad $H$ and $H$ is connected, $H$ leaves each $V_j$ invariant, and the subgroup $M_j$ generated by the subgroups $U_y$ $(y \in V_j)$ is normal in $H$. We claim that $M_i = H$ for some $i \in I$. This is obviously the case for any $i$ if $H$ is simple. The group $H$ is simple except when $n = 4$ and $F$ has signature $(2, 2)$. Then $\mathfrak{h}$ is the direct sum of two copies of the Lie algebra of $\mathbf{SL}_2(\mathbf{R})$. In that case, we may take for $V_i$ the Zariski closure of the orbit of any element $y \in \mathcal{N}$ whose projections on the two factors are not zero.

For $y \in \mathcal{N}$, let $L_y$ be the closed subgroup such that $L_y.o$ is the closure of $U_y.o$. Let $\mathcal{L}$ be the set of the subgroups $L_y$. For $L \in \mathcal{L}$, let $M_L$ be the set of $y \in V_i$ such that $L_y.o \subset L.o$. It is obviously closed. Since $\mathcal{L}$ is countable, there exists at least one $L$ such that $M_L$ contains a non-empty open subset of $V_i$. Choose one such $L$ of smallest possible dimension. Then let $Y$ be a non-empty open subset of $U_i$ such that $\overline{U_y.o} \subset L.o$ for $y \in Y$. Let $R$ be the subgroup generated by the groups $U_y$ $(y \in Y)$. Then any $r \in R$ leaves $L.o$ stable, hence $R \subset L$ and $\overline{R.o} \subset L.o$. Let $\mathfrak{r}$ be the Lie algebra of $R$. Then $\mathfrak{r} \cap V_i$ is an algebraic subset which contains a non-empty open subset of $V_i$, hence is equal to $V_i$. Therefore $R = H$ and $L.o$ is the closure of $H.o$.

**3.7.** *Remarks.* (i) The argument in (a) is valid if $G$ is a linear algebraic group defined over $\mathbf{Q}$ and $\Gamma$ an arithmetic subgroup.

(ii) The proof in (b) is valid without change in the general case, once the existence of $V_i$ is proved. It can be easily deduced from the fact that $H$ is the semi-direct product of a normal nilpotent subgroup, all of whose elements are Ad-unipotent, and of a semisimple group without compact factors. For another argument, see [R4].

**3.8.** The crucial difference between the case $n = 2$, where the Oppenheim conjecture is false, and $n \geqq 3$ lies in the fact that $SO(F)^o$ is generated by unipotent elements for $n \geq 3$ but does not contain any (except 1) for $n = 2$. In fact, in that last case, the flow defined by $H$ is the geodesic flow, and it is well known that its orbits may be neither dense nor closed and may have closures which are not manifolds.

**3.9.** Let $H = \exp \mathbf{R}.y$ with $y \in \mathfrak{g}$ ad-nilpotent be an Ad-unipotent one-parameter group. In this case, M. Ratner establishes a further property of $H.x$, namely, that $H.x$ is *uniformly distributed* in its closure [R4, Theorem B]:

*Let $L = L(x, H)$ be as in the remark to Theorem 3 and $d\nu_L$ the $L$-invariant*

*probability measure with support* $Lx$. *Then*

$$(1) \qquad \lim_{t\to\infty} t^{-1} \int_0^t f\big((\exp s.y).x\big) ds = \int f.d\nu_L$$

*for every bounded continuous function* $f$ *on* $\Omega$.

This had been proved for $G = \mathbf{SL}_2(\mathbf{R})$ in [DS] and for $G$ nilpotent in [P]. An extension of Ratner's result to connected unipotent subgroups has been given by N. Shah [Sh, Corollary 1.3].

## 4. **The measure theoretic conjecture.**

We have so far emphasized the topological conjecture because of its relevance to the proof of the Oppenheim conjecture. However, it appeared comparatively recently in ergodic theory (motivated by the "Davenport conjecture" in fact), in the context of activity centering on a basic problem in ergodic theory: given a group $L$ acting on a measure space $X$, classify the ergodic $L$-invariant probability measures. To complete the picture, I will now discuss one aspect of this problem. From the point of view of the applications to the Oppenheim conjecture, it is not strictly needed, as pointed out in 2.6. However, it is an essential (and the hardest) step in the work of M. Ratner leading to 3.3 (and 3.9).

4.1. The starting point here is again the horocyclic flow on a Riemann surface of finite area (see 3.1; we use the same notation). When $G/\Gamma$ is compact, H. Furstenberg [F] proved, as a strengthening of the fact that all orbits of $H$ are dense, that the horocycle flow is "uniquely ergodic" (only one $H$-invariant ergodic probability measure). The existence of closed orbits when $G/\Gamma$ is not compact shows this is not so in general. But the results of [D1] imply that an $H$-invariant ergodic probability measure on $G/\Gamma$ is either $G$-invariant or supported by a closed orbit of $H$. This led to the following conjecture, to be called here the "measure theoretic conjecture".

4.2. **Conjecture** (Dani [D2], Margulis [M1, M3]). *Let* $H$ *be Ad-unipotent and* $\mu$ *an ergodic* $H$-*invariant probability measure on* $\Omega$. *Then there exists a closed subgroup* $L$ *of* $G$ *containing* $H$, *a point* $x \in \Omega$ *such that* $L.x$ *is closed and* $\mu$ *a* $L$-*invariant measure with support* $L.x$.

To be more precise, this is conjectured in [D2] when $G$ is reductive (as in 3.2) and $H$ one-dimensional. Moreover, it is proved in [D2] when $H$ is a maximal horospherical subgroup of $G$.

In her papers M. Ratner refers to this, or rather to a variant of it (see below) as the "measure theoretic Raghunathan conjecture" because it is the measure theoretic analogue of the topological conjecture. It seems to me this is somewhat of a misnomer, since, as far as I know, Raghunathan did not consider the measure theoretic case at all.

4.3. To state Ratner's theorem, we use a definition introduced in [R1, R2, R3]: let $\mu$ be an $H$-invariant probability measure on $\Omega$. Denote by $\Lambda = \Lambda(\mu)$ the set of $g \in G$ which leaves it invariant. It is a closed subgroup [Rl, Proposition 1.1]. Then $\mu$ is said to be *algebraic* if there exists $x \in \Omega$ such that $G_x \cap \Lambda$ is of finite covolume in $\Lambda$ and $\Lambda x$ is the support of $\mu$. In particular $\Lambda x$ is closed [R, Theorem 1.13].

**Theorem 4** (M. Ratner). *We drop the assumption that $\Gamma$ has finite covolume. Let $H$ be as in Theorem 3. Then any ergodic $H$-invariant probability measure $\mu$ on $\Omega$ is algebraic. If $H$ is connected, it contains a one-parameter subgroup, Ad-unipotent in $G$, which acts ergodically on $(\Omega, \mu)$.*

The first assertion is proved in [R3, Theorem 3], after having been established for $G$ solvable in [R1] and for $G$ semisimple, $\Gamma$ cocompact, in [R2]. The second one is Proposition 5.2 in [R3].

To reduce the proof to the one-dimensional case, [R3] also provides a countability statement (Theorem 2 there), namely:

**Theorem 5.** *Fix $x \in \Omega$. Let $\Phi_x(G, \Gamma)$ be the set of closed connected subgroups $L$ of $G$ with the following property: $G_x \cap L$ has finite covolume in $L$, and $L$ contains a connected subgroup of $G$ generated by Ad-unipotent elements of $G$ which acts ergodically on $(L.x, \nu_L)$, where $\nu_L$ is the $L$-invariant probability measure on $L.x$. Then $\Phi_x(G, \Gamma)$ is countable.*

## 5. Some remarks on the proofs and further developments.

5.1. We shall not try to describe the proofs, which take over 200 pages, and limit ourselves to some comments, all the more since we can refer to Ratner's survey [R9] for further information.

There is so far only one proof of Theorem 4, given in [R1, R2, R3]. It is also described for $\mathbf{SL}_2(\mathbf{R})$ in [R5] and sketched in the general case in [R6]. The assumption Ad-unipotent for one-parameter subgroup $H$ is used in two crucial ways. First, the adjoint representation Ad $\mathfrak{g}: h \mapsto \mathrm{Ad}_{\mathfrak{g}} h$ of $H$ on the Lie algebra $\mathfrak{g}$ of $G$ is given by a polynomial mapping of $H$ into $\mathrm{End}(\mathfrak{g})$, and all orbits of $H$ there are closed. Second, if $G$ is semisimple, the Lie algebra $\mathfrak{h}$ of $H$ belongs to a " $\mathfrak{sl}_2$-triple"; i.e. there exists a homomorphism $\varphi: M \to G$ of a covering $M$ of $\mathbf{SL}_2(\mathbf{R})$ such that $H$ is the (isomorphic) image of the identity component of the inverse image of the group of upper triangular unipotent matrices of $\mathbf{SL}_2(\mathbf{R})$. The image $A$ under $\varphi$ of the identity component of the inverse image of the group of diagonal matrices in $\mathbf{SL}_2(\mathbf{R})$ then normalizes $H$, and Ad $\mathfrak{g} A$ is diagonalisable (over the reals). M. Ratner then says that $A$ is diagonal, or is a diagonal, for $H$. This allows one in particular to use the representation theory of $\mathbf{SL}_2(\mathbf{R})$ to describe the actions of $A$ and $H$ on $\mathfrak{g}$ by the adjoint representation.

The first fact yields some control of some orbits, in particular of the time passed in certain subsets. The starting point of such estimates is the following property of polynomials on the line:

Let $\mathcal{P}(n)$ be the set of polynomials on $\mathbf{R}$ of degrees $\leqq n$. Then there exists $\eta \in (0, 1)$ with the following property: if $P \in \mathcal{P}(n)$ is such that for given $t$, $\theta > 0$, it satisfies the condition

$$\max_{s \in [0, t]} |P(s)| = |P(t)| = \theta,$$

then $\theta/2 < |P(s)| \leq \theta$ for $s \in [(1 - \eta)t, t]$.

5.2. The passage from Theorem 4 to the uniform distribution theorem 3.9 is carried out in [R4], described for $\mathbf{SL}_2(\mathbf{R})$ in [R5] and sketched for the general case in [R6]. The proof is by induction on $\dim G$ so that it may be assumed that there is no proper closed connected subgroup $M$ of $G$ containing $H$ such that $M \cap G_x$ is of finite covolume in $M$. In this case, it has to be shown that $H.x$ is uniformly distributed with respect to the $G$-invariant probability measure $d\nu_G$ on $\Omega$. This will also imply that $H.x$ is dense in $G$. For $t > 0$, let $T_{x,\,t}$ be the measure on the space $C_o(G)$ of bounded continuous functions on $G$ defined by

$$T_{x,\,t} = t^{-1} \int_0^t f\big(\exp sy\big).x\big)ds.$$

Then it has to be shown that

$$d\nu_G = \lim_{t \to \infty} T_{x,\,t}$$

in the weak sense, i.e. $\int f d\nu_G = \lim T_{x,\,t}(f)$ for $f \in C_o(G)$. The measures $T_{x,\,t}$ are obviously $\leqq 1$ in norm. Since a bounded set of measures in the weak * topology is relatively compact, the set $M(x,\,H)$ of measures which are limit points of sequences $T_{x,\,t_j}$ $(t_j \to \infty)$ is not empty. All these measures are $H$-invariant. One has to prove eventually that $M(x,\,H)$ consists solely of $d\nu_G$. Let $\mu \in M(x,\,H)$ and $Y$ be its support. It is shown first that $\mu(\Omega) = 1$. By a general fact, $\mu$ admits a decomposition into ergodic $H$-invariant probability measures, i.e. there exists a family of $\xi = \{\mu_y\}$ of ergodic $H$-invariant probability measure so that the supports $C(y)$ of the $\mu_y$ form a partition of $Y$ and a measure $\nu_\xi$ on the quotient $Y/\xi$ such that $\mu(f) = \int \mu_y\big(f|C(y)\big)\nu_\xi$. (All this is not really true as stated, but only up to sets of measure zero for the measures under consideration.) Now by Theorem 4, each $\mu_y$ is algebraic, i.e. there exists for $y \in Y$ a closed subgroup $\Lambda_y$ containing $H$ such that $C(y) = \Lambda_y^o.y$, the intersection $\Lambda_y^o \cap G_y$ has finite covolume in $\Lambda_y^o$ and $\mu_y$ is the $\Lambda_y^o$-invariant probability measure on $\Lambda_y^o.y$. Then the main part of the proof consists in showing, under our initial assumption, that the $\mu$-measure of the union of the $C(y)$ which are $\neq \Omega$ is zero. This is established for $G$ semisimple in Theorem 2.1 and in the general case in Corollary 3.1 of [R4].

5.3. Another way to go from Theorem 4 to 3.9 is described in [DM4]. The authors also prove their own variant of a countability theorem (Theorem 5.1):

*Fix a right invariant Riemannian metric on $G$, whence also a Riemannian metric on $\Omega$. Given $c > 0$, let $\mathcal{V}_c$ be the set of closed connected subgroups $H$ such that $H\Gamma/\Gamma$ is closed and has volume $\leqq c$. Then the set of intersections $H \cap \Gamma$ for $H \in \mathcal{V}_c$ is finite. Let further $\rho : G \to GL(V)$ be a finite dimensional representation of $G$ with kernel central in $G$. Then the set of $H \in \mathcal{V}_c$ for which $\rho(H \cap \Gamma)$ is Zariski dense in $\rho(H)$ is finite.*

In this theorem, $\Gamma$ need not have finite covolume.

5.4. In fact, [DM4] proves a generalization of 3.9, also independently obtained in [R6, Theorem 7], involving a sequence of Ad-unipotent subgroups. Let $U_n$ $(n \in \mathbf{N})$ and $U$ be Ad-unipotent one-parameter subgroups of $G$. The relation $U_n \to U$ means, by definition, that $U_n(t) \to U(t)$ for every $t \in \mathbf{R}$. We refer to 3.3 for the definition of $L(x,\,U)$. In view of Theorem 3, $L(x,\,U) = G$ if and only if $U.x$ is dense in $\Omega$.

**Theorem 6.** *Let $U$ be a one-dimensional Ad-unipotent subgroup of $G$ such that $L(x,\,U) = G$ and $x \in \Omega$. Let $x_n \in \Omega$ tend to $x$ and $U_n$ be a sequence of one-dimensional Ad-unipotent subgroups of $G$ tending to $U$. Then for any sequence $t_n \to \infty$ and any bounded continuous function $f$ on $\Omega$ we have*

$$\lim_{t_n \to \infty} t_n^{-1} \int_0^{t_n} f(U_n(s).x_n).ds = \int_\Omega f d\nu_G.$$

5.5. In [MS] the authors consider a sequence of Ad-unipotent one-dimensional subgroups $\{U_n\}$ (not necessarily convergent) and a convergent sequence of measures $\mu_n \to \mu$, where $\mu_n$ is ergodic and $U_n$-invariant. They show, among other things, that $\mu$ is algebraic, invariant and ergodic for some one-dimensional Ad-unipotent subgroup.

5.6. The assumption that $H^o$ is generated by Ad-unipotent subgroups is of course essentially used. Some assumption on $H$ is certainly needed since, for instance, some orbits of the group of diagonal matrices $A$ of $\mathbf{SL}_2(\mathbf{R})$ in $\mathbf{SL}_2(\mathbf{R})/\mathbf{SL}_2(\mathbf{Z})$ have closures which are not even manifolds. Nevertheless, Theorems 3 and 4 (resp. Theorem 4) have been extended to some more general classes of groups by M. Ratner (resp. S. Mozes).

a) [R6, Theorem 9]. $H$ is connected, generated by a closed connected subgroup $M$, itself generated by Ad-unipotent one-dimensional subgroups, and by subgroups $A_i$ $(1, \ldots, m)$ where $A_i$ is diagonal with respect to some one-dimensional Ad-unipotent subgroup $U_i$ of $M$ (see 5.1 for that notion).

b) [Mo]. $G$ has a connected semi-simple subgroup $L$ without compact factors containing $H$, and $H$ is connected, epimorphic in $L$.

Here epimorphic is in the sense of [BB]. This is equivalent to requiring that the regular functions on $L/H$ be only the constants.

There is an overlap between these two classes. For instance, both contain the parabolic subgroups of a connected semisimple subgroup $L$ of $G$ without compact factors.

### III. The S-arithmetic case

The Oppenheim conjecture gives a criterion for an indefinite quadratic form to be "rational", meaning rational with respect to $\mathbf{Q}$. In [RR], the authors initiated the consideration of an analogous question over a number field $k$. It involved looking at quadratic forms over the archimedean completions of $k$. This was taken up and generalized in [BP1] and [BP2], where finite places are also included (following a suggestion of G. Faltings). Subsequently, extensions of some of (resp. all) the results on flows have been obtained in [MT2] (resp. [R8]). To complete the picture, I describe some of these generalizations in this section, assuming familiarity with some basic concepts in algebraic number theory and also, in §§7 and 8, with the theory of linear algebraic groups.

### 6. The generalized Oppenheim conjecture.

6.1. In the sequel, $k$ is a number field and $\mathfrak{o}$ the ring of integers of $k$. For every normalized absolute value $|\cdot|_v$ on $k$, let $k_v$ be the completion of $k$ at $v$. In the sequel $S$ is a finite set of places of $k$ containing the set $S_\infty$ of the archimedean ones, $k_S$ the direct sum of the fields $k_s$ $(s \in S)$ and $\mathfrak{o}_S$ the ring of $S$-integers of $k$ (i.e. of elements $x \in k$ such that $|x|_v \leqq 1$ for $v \notin S$). For $s$ non-archimedean, the valuation ring of $k_s$ is denoted $\mathfrak{o}_s$.

Let $F$ be a quadratic form on $k_S^n$. Equivalently, $F$ can be viewed as a family $(F_s)$ $(s \in S)$, where $F_s$ is a quadratic form on $k_s^n$. The form $F$ is non-degenerate if and only each $F_s$ is non-degenerate. We say that $F$ is isotropic if *each* $F_s$ is so, i.e. if there exists for every $s \in S$ an element $x_s \in k_s^n - \{0\}$ such that $F_s(x_s) = 0$.

The form $F$ is said to be *rational* (over $k$) if there exists a quadratic form $F_o$ on $k^n$ and a unit $c$ of $k_S$ such that $F = c.F_o$, *irrational* otherwise.

6.2. The following theorem reduces to the truth of the Oppenheim conjecture if $k = \mathbf{Q}$ and $S = S_\infty$ consists of the infinite place. That (i) and (ii) are equivalent is the generalized Oppenheim conjecture.

**Theorem 7** ([BP2], Theorem A). *Let $n \geqq 3$ and $F$ be an isotropic non-degenerate quadratic form on $k_S^n$. Then the following two conditions are equivalent*:

(i) *$F$ is irrational.*

(ii) *Given $\epsilon > 0$, there exists $x \in \mathfrak{o}_S^n$ such that $0 < |F_s(x)|_s < \epsilon$ for all $s \in S$.*

Here again (ii) $\Rightarrow$ (i) is obvious, and the main point is (i) $\Rightarrow$ (ii). For $S = S_\infty$, the proof is patterned after that of Margulis [M3]; it is based on a generalization of Proposition 1 (and of Theorem 1$'$ in [M3]) to the case where

$$(1) \qquad G = \prod_{s \in S} \mathbf{SL}_3(k_s), \ \Gamma = \mathbf{SL}_3(\mathfrak{o}_S), \ H = \prod_{s \in S} SO(F_s).$$

To treat the general case, this result is combined with an argument using strong approximation in algebraic groups and some geometry of numbers.

*Remarks*. 1) We have assumed that each $F_s$ is isotropic over $k_s$. If this is not so, it is easily seen that Theorem 7 cannot hold (see 1.10 in [BP2]).

2) In the original case $k = \mathbf{Q}$ and $S = S_\infty$, we already pointed out in 1.3 that the truth of the Oppenheim conjecture and one theorem of Oppenheim imply that $F(\mathbf{Z}^n)$ is dense in $\mathbf{R}$ when $F$ is irrational. In the general case, going from Theorem 7 to the density appears to be more difficult. At this time, it has to make use of the orbit closure theorem, which then yields again a much stronger statement (see §8).

## 7. **Closures of orbits.**

7.1. The results of §6 and their proofs led one naturally to ask whether the results on flows reviewed in II would extend to a framework including the one of 6.2, where $G$ would be a product of real and $p$-adic groups. This generalization was carried out by M. Ratner [R7, R8, R9] for all of her results and, independently, by G.A. Margulis and G. Tomanov [MT1, MT2] for the measure theoretic conjecture in the setting of algebraic groups. We focus here on the orbit closure theorem and then, in the next section, deduce from it an extension of Proposition 2, hence also of Theorem 2, to the $S$-arithmetic case.

7.2. Again, the framework of Ratner's work is Lie group theory, rather than algebraic groups. We refer to [Bu1] or [S1] for the notion of Lie group over a local field $E$ of characteristic zero (i.e. $\mathbf{R}$, $\mathbf{C}$ or a finite extension of $\mathbf{Q}_p$ or, equivalently, the fields $k_s$ of 6.2, for variable number fields $k$). If $\mathcal{G}$ is an algebraic group defined over $E$, then the group $\mathcal{G}(E)$ of rational points of $\mathcal{G}$ is in a natural way a Lie group over $E$, and the Lie algebra of $\mathcal{G}(E)$, as a Lie group, is the space of rational points over $E$ of the Lie algebra of $\mathcal{G}$, as an algebraic group.

7.3. Let $k$, $S$, $k_s$, $o_s$, $k_S$ be as in 6.1. For each $s \in S$, there is given a Lie group $G_s$ over $k_s$ and a closed subgroup $H_s$ generated by Ad-unipotent one-dimensional subgroups over $k_s$. The product $G$ of the $G_s$ is then in a natural way a locally compact topological group, and the product $H$ of the $H_s$ is a closed subgroup.

As usual, we identify $G_t$ $(t \in S)$ to the subgroup of $G$ consisting of the elements $(g_s)_{s \in S}$ such that $g_s = 1$ for $s \neq t$.

Two slight restrictions are imposed on $G_s$ if $s$ is non-archimedean. First, the kernel of the adjoint representation is the center $Z(G_s)$ of $G$. (If $G_s$ is the group of rational points of an algebraic group which is connected in the Zariski topology, this is automatic, otherwise the kernel could be bigger.) Second, it is required that the orders of the finite subgroups of $G_s$ are bounded. This is always true if $G_s$ is linear. An argument is given in [S1, IV, Appendix 3], proof of Theorem 1. We sketch it: the maximal compact subgroup of $\mathbf{GL}_n(k_s)$ are all conjugate to $\mathbf{GL}_n(o_s)$, [S1, Theorem 1, p. 122], so it suffices to consider the finite subgroup of $\mathbf{GL}_n(o_s)$. Since $\mathbf{GL}_n(o_s)$ is compact, it follows from [S1, Theorem 5, p. 119] that it contains a torsion-free normal open subgroup $N$. Then $\mathbf{GL}_n(o_s)/N$ is a finite group, and any finite subgroup of $\mathbf{GL}_n(o_s)$ is isomorphic to a subgroup of that quotient.

In [R8], $G_s$ is said to be Ad-regular if it satisfies the first condition and regular if it satisfies both. In particular, we see from the above that if $G_s$ is the group of rational points of a connected linear algebraic group defined over $k_s$, it is regular.

In Theorem 8 below, $G_s$ is assumed to be regular for $s \in S$ non-archimedean.

**Theorem 8** [R8, Theorem 2]. *Let $G$, $H$ be as above. Let $M$ be a closed subgroup of $G$ containing $H$ and $\Gamma$ a discrete subgroup of finite covolume of $M$. Let $x \in M/\Gamma$. Then $M$ contains a closed subgroup $L$ such that $L.x$ is the closure of $H.x$ and $L \cap M_x$ has finite covolume in $L$.*

This is stated and proved by M. Ratner for $k = \mathbf{Q}$, but this is no loss in generality. Let $\mathbf{Q}_s$ be the completion of $\mathbf{Q}$ in $k_s$. It is therefore equal to $\mathbf{R}$ if $k_s = \mathbf{R}$, $\mathbf{C}$ and to a field $\mathbf{Q}_p$ for some prime $p$ if $k_s$ is non-archimedean. Then a Lie group over $k_s$, of dimension $m$, may be viewed in a natural way as a Lie group over $\mathbf{Q}_s$, of dimension $m[k_s : \mathbf{Q}_s]$, in the same way as a complex Lie group can be viewed as a real Lie group of twice the dimension. Let us denote by $G'_s$ the group $G_s$ thus endowed with a structure of Lie group over $\mathbf{Q}_s$. Then the identity map of the product $G'$ of the $G'_s$ onto $G$ is an isomorphism of topological groups. Moreover, if $U$ is a one-dimensional Ad-unipotent group over $k_s$, then $U'$ is a direct sum of $[k_s : \mathbf{Q}_s]$ one-dimensional Ad-unipotent subgroup over $\mathbf{Q}_s$ of $G'_s$. Therefore it is clear that Theorem 8 for $G'$ implies it for $G$.

One advantage of the shift to Lie groups over $\mathbf{Q}_s$ is the fact (due to E. Cartan over the real numbers) that any closed subgroup of a Lie group over $\mathbf{Q}_s$ is a Lie group over $\mathbf{Q}_s$ [Bu1, Chapter 3, §8 no. 2; S1, Chapter V, §9].

7.4. The steps in the proof of Theorem 8 are similar to those for real Lie groups. There is first a theorem proving that $H$-invariant ergodic probability measures are algebraic, also established in [MT2] when the $G_s$ are groups of rational points of linear algebraic groups. Then a countability statement [R8, Theorem 1.3] allows one to reduce the proof of Theorem 8 to the case where $H$ is one-dimensional Ad-unipotent, contained in one factor, in which case a uniform distribution theorem is also proved [R8, Theorem 3].

8. **Applications to quadratic forms.** We now generalize the density theorems of 2.3 and 2.5 to the present case, using Theorem 8, in the same way as was done in the case $S = S_\infty$ in §7 of [BP2]. As stated in the introduction, we include the proofs, obtained jointly with G. Prasad.

8.0. We shall use the following lemma. It should be known in the theory of affine symmetric spaces. For lack of a reference, we have included a proof in the appendix.

**Lemma.** *Let $E$ be a field of characteristic zero, $\mathfrak{g}$ a simple Lie algebra over $E$, $\sigma \neq 1$ an involutive automorphism of $\mathfrak{g}$ and $\mathfrak{k}$ the fixed point set of $\sigma$. Assume that $\mathfrak{k}$ is semi-simple. Then any $\mathfrak{k}$-invariant subspace of $\mathfrak{g}$ containing $\mathfrak{k}$ is equal to $\mathfrak{k}$ or $\mathfrak{g}$. In particular, $\mathfrak{k}$ is a maximal proper subalgebra of $\mathfrak{g}$.*

8.1. We revert to the notation of 6.2. Moreover, let

$$G_s = \mathbf{SL}_n(k_s), \ H_s = SO(F_s) \qquad G = \prod_{s \in S} G_s, \qquad H = \prod_{s \in S} H_s,$$

$\mathcal{G}_s$ be $\mathbf{SL}_n$ viewed as algebraic group over $k_s$ and $\mathcal{H}_s$ the algebraic group over $k_s$ such that $\mathcal{H}_s(k_s) = H_s$.

Following a notation of [BT], we let $H_s^+$ denote the subgroup of $H_s$ generated by one-dimensional unipotent (hence Ad-unipotent) subgroups. We claim that it is a closed and open normal subgroup of finite index of $H_s$. If $k_s = \mathbf{C}$, this is immediate, since $H_s$ is semisimple and connected in the usual topology, and in fact $H_s = H_s^+$. If $k_s = \mathbf{R}$, then $H_s^+$ is the topological identity component of $H_s$ and has index two. Now let $k_s$ be non-archimedean. Let $\tilde{\mathcal{H}}_s$ be the universal covering of $\mathcal{H}_s$, i.e. the spinor group of $F_s$, and $\mu\colon \tilde{\mathcal{H}}_s \to \mathcal{H}_s$ the central isogeny. Let $\tilde{H}_s = \tilde{\mathcal{H}}_s(k_s)$. It is known that $\tilde{H}_s = \tilde{H}_s^+$ is generated by one-dimensional unipotent subgroups [BT, 6.15], that $\mu(\tilde{H}_s^+) = H_s^+$ [BT, 6.3] and that $\mu(\tilde{H}_s)$ is a normal open and closed subgroup of finite index of $H_s$ [BT, 3.20], whence our assertion in that case.

We note that $H_s^+$ is not compact, since it is of finite index in $H_s$ and the latter, being the orthogonal group of an *isotropic* form, is not compact.

We let $\mathfrak{h}_s$ be the Lie algebra of $H_s$ and $N_s$ the normalizer of $\mathfrak{h}_s$ in $G_s$, i.e.

$$N_s = \{g \in G_s | \operatorname{Ad} g(\mathfrak{h}_s) = \mathfrak{h}_s\}.$$

We claim that $N_s$ is also the normalizer of $H_s$ or of $H_s^+$. In fact, both groups, viewed as Lie subgroups of $G_s$, have $\mathfrak{h}_s$ as their Lie algebra, therefore any element $g \in G_s$ normalizing $H_s$ or $H_s^+$ belongs to $N_s$. Conversely, let $g \in N_s$. Since $\mathfrak{h}_s$ is the space of rational points of the Lie algebra of $\mathcal{H}_s$ and is of course Zariski-dense in it, the automorphism $\operatorname{Int} g\colon x \mapsto g.x.g^{-1}$ of $\mathcal{G}_s$ leaves $\mathcal{H}_s$ stable, hence $g$ normalizes $H_s$ and therefore also $H_s^+$.

**Lemma.** (i) $H_s^+$ *has finite index in $N_s$.* (ii) *Let $M$ be a subgroup of $G_s$ containing $H_s^+$. Then either $M = G_s$ or $M \subset N_s$.*

(i) Since $H_s^+$ has finite index in $H_s$, it suffices to show that $H_s$ has finite index in $N_s$. The only quadratic forms on $k_s^n$ invariant under $H_s$ are the multiples of $F_s$. If $x \in N_s$, then ${}^t x.F_s.x$ is invariant under $H_s$, hence of the form $c.F_s$ ($c \in k_s^*$). It has the same determinant as $F_s$; hence $c^n = 1$, and therefore $N_s/H_s$ is isomorphic to a subgroup of the group of $n$-th roots of unity.

(ii) Identify $F_s$ to a symmetric, invertible, matrix. Then the map

$$\sigma\colon x \mapsto F_s.{}^t x^{-1}.F_s^{-1} \qquad (x \in G_s)$$

is an automorphism of $G_s$, obviously of order two, and $H_s$ is the fixed point set of $\sigma$. The differential $d\sigma$ of $\sigma$ at the origin is an involutive automorphism of $\mathfrak{g}_s$ with fixed point set $\mathfrak{h}_s$. The group $\mathcal{G}_s$ (resp. $\mathcal{H}_s$) is simple (resp. semisimple) as an algebraic group; therefore $\mathfrak{g}_s$ (resp. $\mathfrak{h}_s$) is a simple (resp. semisimple) Lie algebra. By 8.0, any $\mathfrak{h}_s$-invariant subspace of $\mathfrak{g}_s$ containing $\mathfrak{h}_s$ is equal to $\mathfrak{h}_s$ or to $\mathfrak{g}_s$.

Now let $M$ be a subgroup of $G_s$ containing $H_s^+$ but not contained in $N_s$. We have to show that $M = G_s$. Let $\mathfrak{g}$ be the subspace generated by the subalgebras $\operatorname{Ad} m(\mathfrak{h}_s)$, $(m \in M)$. It is normalized by $M$, obviously, and in particular by $H_s^+$. Therefore it is $\mathfrak{h}_s$-invariant. It is $\neq \mathfrak{h}_s$, since $M$ is not in $N_s$. By the remark just made, $\mathfrak{g} = \mathfrak{g}_s$. There exists therefore a finite set of elements $m_i \in M \, (1 \leq i \leq a)$ such that $\mathfrak{g}_s$ is the sum of the subalgebras $\mathfrak{h}_i = \operatorname{Ad} m_i(\mathfrak{h}_s)$. The Lie algebra $\mathfrak{h}_i$ is the Lie algebra of $H_i = m_i.H_s^+.m_i^{-1}$. Let $Q = H_1 \times \ldots \times H_d$ be the product of the $H_i$ and $\mu : Q \to G_s$ be the map which assigns to $(h_1, \ldots, h_a)$ $(h_i \in H_i)$ the product of the $h_i$'s. It is a morphism of $k_s$-manifolds, whose image is contained in $M$. The tangent space at the identity of $Q$ is the direct sum of the $\mathfrak{h}_i$'s. Therefore the differential $d\mu$ of $\mu$ at the identity maps the tangent space to $Q$ onto $\mathfrak{g}_s$. This implies that $\mu(Q)$ contains an open neighborhood of the identity in $G_s$ (see [S1, III, 10.2]). Since it belongs to $M$, the latter is an open subgroup of $G_s$. It contains $H_s^+$, which is not compact, as noted in 8.1, hence is noncompact. Moreover, it is elementary that $G_s = \mathbf{SL}_n k_s$ is generated by the group of unipotent upper triangular matrices and its conjugates. It then follows from Theorem (T) in [Pr] that $M = G_s$.

8.3. Let $\Gamma = \mathbf{SL}_n(\mathfrak{o}_S)$. It is viewed as a discrete subgroup of $G$ via the embeddings $\mathbf{SL}_n(k) \to \mathbf{SL}_n(k_s)$. The quotient $\Omega = G/\Gamma$ has finite volume. We let $o$ be the coset $\Gamma$ in $\Omega$.

**Theorem 9.** *If $F$ is irrational, the orbit $H.o$ is dense in $\Omega$.*

Let $H^+$ be the product of the groups $H_s^+$ (see 8.1). Since $H_s^+$ has finite index and is normal, open and closed in $H_s$, the same is true for $H^+$ in $H$, and it is equivalent to prove that $H^+.o$ is dense in $\Omega$.

By Theorem 8 (with $M = G$), there exists a closed subgroup $L$ of $G$ such that $L.o$ is the closure of $H^+.o$ and $L \cap \Gamma$ has finite covolume in $L$.

Let $M_s = L \cap G_s$. It is a closed normal subgroup of $L$ which contains $H_s^+$. By 8.2, we have either $M_s \subset N_s$ or $M_s = G_s$. Now let $P_s$ be the projection of $L$ into $G_s$. It normalizes $M_s$ and contains it. Assume $M_s \subset N_s$. Then $\operatorname{Ad} g \, (g \in P_s)$ leaves invariant the Lie algebra of $M_s$, which is the same as that of $H_s$, hence $g$ belongs to $N_s$. In particular, $P_s$ is closed and open in $N_s$. If $M_s = G_s$, then $P_s = G_s$. Therefore the product $M$ of the $M_s$ is normal, closed and open, of finite index in the product $P$ of the $P_s$, and $P$ is closed. We have of course $M \subset L \subset P$. As a consequence, $M$ is normal, open and closed, of finite index, in $L$.

Now define $Q_s$ by the rule: $Q_s = H_s$ if $M_s \subset N_s$, and $Q_s = G_s$ if $M_s = G_s$, and let $Q$ be the product of the $Q_s$. Then $Q \cap L$ is open and closed, of finite index, in both $L$ and $Q$. Therefore $Q \cap \Gamma$ has finite covolume in $Q$. By Proposition 1.2 in [BP2], there exists a $k$-subgroup $\mathcal{Q}$ of $\mathbf{SL}_n$ such that $\mathcal{Q}(k_s) = Q_s$ for every $s \in S$. This shows first of all that either $Q_s = G_s$ for all $s \in S$ or $Q_s = H_s$ for all $s \in S$. In the first case, $L = G$ and $H^+.o$ is dense. We have to rule out the second one. In that case $H.o$ is closed, $L = H^+$ and $H \cap \Gamma$ has finite covolume in $H$. Moreover $\mathcal{Q}$ is the orthogonal group of a form $F_o$ on $k^n$, and there exists a unit $c$ of $k_S$ such that $F = c.F_o$, i.e. $F$ is rational over $k$, contradicting our assumption.

8.4. We can now generalize 2.5.

A subset $(x_1, \dots, x_m)$ $(m \leqq n)$ of $\mathfrak{o}_S^n$ is said to be primitive if it is part of a basis of $\mathfrak{o}_S^n$ over $\mathfrak{o}_S$. If $m < n$, it is so if and only if there exists $g \in \mathbf{SL}_n(\mathfrak{o}_S)$ such that $x_i = g(e_i)$ $(i = 1, \dots, m)$, where $e_i$ is the $i$-th canonical basis element of $k_S^n$.

**Corollary.** *Assume $F$ to be irrational. Let $\lambda_i \in k_s$ $(i = 1, \dots, n-1)$. Then there exists a sequence of primitive $(n-1)$-tuples $(x_{j,\,1}, \dots, x_{j,\,n-1})$ $(j = 1, 2 \dots)$ in $\mathfrak{o}_S^n$ such that*

$$\lambda_i = \lim_{j \to \infty} F(x_{j,\,i}) \qquad (i = 1, \dots, n-1).$$

*In particular, the set of values of $F$ on primitive elements of $\mathfrak{o}_S^n$ is dense in $k_S$.*

The argument is the same as in 7.9 in [BP2]. We repeat it for the sake of completeness.

Let $\lambda_{i,\,s}$ be the component of $\lambda_i$ in $k_s$ $(s \in S)$. The form $F_s$, being isotropic, takes all values in $k_s$. The representation of $H_s$ in $k_s^n$ is irreducible, and no level surface $F_s = c$ is contained in a hyperplane; hence, given $s \in S$, we can find linearly independent vectors $y_{s,\,i} \in k_s^n$ such that $F_s(y_{s,\,i}) = \lambda_{i,\,s}$. There exists then $g_s \in G_s$ such that $g_s(e_i) = y_{s,\,i}$ $(i = 1, \dots, n-1)$. Let $g = (g_s)$ and $y_i = (y_{s,\,i}) \in k_S^n$. Then

$$F(y_i) = F(g.e_i) = \lambda_i \qquad (i = 1, \dots, n-1).$$

By Theorem 8, $H.o$ is dense in $\Omega$, hence $H.\Gamma$ is dense in $G$. There exist therefore elements $h_j \in H$, $\gamma_j \in \Gamma$ $(j = 1, 2, \dots)$ such that $h_j.\gamma_j \to g$. Then

$$\lambda_i = F(g.e_i) = \lim_{j \to \infty} F(h_j.\gamma_j.e_i) = \lim_{j \to \infty} F(\gamma_j.e_i) = \lim_{j \to \infty} F(x_{j,\,i}),$$

where $x_{j,\,i} = \gamma_j.e_i$ $(j = 1, 2, \dots; i = 1, \dots, n-1)$. For each $j$, $(x_{j,\,1}, \dots, x_{j,\,n-1})$ is a primitive subset of $\mathfrak{o}_s^n$, whence the corollary.

8.5. **Errata to [BP2].** In the proof of 7.4, the difference between $H^+$ and $H$ has been overlooked. In the archimedean case the subgroup of $H_F$ generated by unipotent elements is the topological identity component $H_F^o$ of $H_F$, and its index is twice the number of real places. The corrections on page 369 are:

Line 2: After $H_s$ add: Let $H_s^o$ be the connected component of the identity in $H_s$ and $H_F^o$ the product of the $H_s^o$.

Line 3: Replace $H_F$ by $H_F^o$.

Line 5: Replace $H_s$ by $H_s^o$.

Line 12: Replace $G$ by $\mathcal{G}$ and 7.2 by 7.3.

Line 23: Replace 7.3 by 7.4.

9. **Appendix: A Lemma on symmetric Lie algebras.** In this appendix, we prove the lemma in 8.0. (We recall that, if $E = \mathbf{R}$ and $\mathfrak{k}$ is maximal compact or $\mathfrak{g}$ is compact, this is true without restriction on $\mathfrak{k}$ and is due to E. Cartan.)

The case where $\mathfrak{g}$ is one-dimensional is left to the reader, so we assume $\mathfrak{g}$ to be also semisimple. Let $\mathfrak{p}$ be the $(-1)$-eigenspace of $\sigma$ in $\mathfrak{g}$. Then

$$(1) \qquad \mathfrak{g} = \mathfrak{k} \oplus \mathfrak{p} \qquad [\mathfrak{k}, \, \mathfrak{k}] \subset \mathfrak{k}, \; [\mathfrak{k}, \, \mathfrak{p}] \subset \mathfrak{p} \qquad [\mathfrak{p}, \, \mathfrak{p}] \subset \mathfrak{k}$$

as usual. If $\mathfrak{m}$ is a $\mathfrak{k}$-invariant subspace of $\mathfrak{p}$, then $\mathfrak{k} \oplus \mathfrak{m}$ is a subalgebra, so that in fact the last clause of the lemma is equivalent to the lemma itself.

Let $B$ be the Killing form of $\mathfrak{g}$. It is non-degenerate and $B(\mathfrak{k}, \, \mathfrak{p}) = 0$; hence the restrictions of $B$ to $\mathfrak{k}$ and to $\mathfrak{p}$ are non-degenerate. We recall that, by invariance

$$(2) \qquad B([a, \, b], \, c) = B(a, \, [b, \, c]) \qquad (a, \, b, \, c \in \mathfrak{g}).$$

The proof is divided into a number of steps:

a) Let $\mathfrak{m} \subset \mathfrak{p}$ be a $\mathfrak{k}$-invariant subspace. Then $[\mathfrak{m}, \, \mathfrak{m}]$ is an ideal of $\mathfrak{k}$ and $[\mathfrak{m}, \, \mathfrak{m}] \oplus \mathfrak{m}$ a subalgebra normalized by $\mathfrak{k}$.

This follows by straightforward application of the Jacobi identity and (1).

b) We have $[\mathfrak{p}, \, \mathfrak{p}] = \mathfrak{k}$.

In fact, $[\mathfrak{p}, \, \mathfrak{p}] \oplus \mathfrak{p}$ is a subalgebra normalized by $\mathfrak{k}$ in view of a), hence a non-zero ideal, hence equal to $\mathfrak{g}$.

c) Let $\mathfrak{m}, \, \mathfrak{n}$ be $\mathfrak{k}$-invariant subspaces of $\mathfrak{p}$ and assume $B(\mathfrak{m}, \, \mathfrak{n}) = 0$. Then $[\mathfrak{m}, \, \mathfrak{n}] = 0$.

Let $m \in \mathfrak{m}, \, n \in \mathfrak{n}$ and $k \in \mathfrak{k}$. Then

$$B(k, \, [m, \, n]) = B([k, \, m], \, n) \subset B(\mathfrak{m}, \, \mathfrak{n}) = 0;$$

therefore $[m, \, n]$ belongs to the radical of the restriction of $B$ to $\mathfrak{k}$. Since the latter is non-degenerate, this proves $[m, \, n] = 0$.

d) There is no proper $\mathfrak{k}$-invariant subspace $\mathfrak{m}$ of $\mathfrak{p}$, on which the restriction of $B$ is non-degenerate.

Let $\mathfrak{m}$ be one and $\mathfrak{n}$ its orthogonal complement. Then $\mathfrak{p} = \mathfrak{m} \oplus \mathfrak{n}$. By c)

$$(3) \qquad [\mathfrak{m}, \, \mathfrak{n}] = 0.$$

By b), $\mathfrak{a} = [\mathfrak{m}, \, \mathfrak{m}] \oplus [\mathfrak{m}]$ and $\mathfrak{b} = [\mathfrak{n}, \, \mathfrak{n}] + \mathfrak{n}$ are subalgebras, and (3) implies that $[\mathfrak{a}, \, \mathfrak{b}] = 0$. By b) and (3), $\mathfrak{k} = [\mathfrak{m}, \, \mathfrak{m}] + [\mathfrak{n}, \, \mathfrak{n}]$; hence $\mathfrak{g} = \mathfrak{a} + \mathfrak{b}$ and $\mathfrak{a}, \, \mathfrak{b}$ are distinct non-zero ideals of $\mathfrak{g}$, a contradiction.

e) There is no proper $\mathfrak{k}$-invariant subspace $\mathfrak{m}$ of $\mathfrak{p}$ on which the restriction of $B$ is degenerate but non-zero.

Assume $\mathfrak{m}$ is such a subspace. Let $\mathfrak{r}$ be the radical of $B|\mathfrak{m}$. It is non-zero, $\neq \mathfrak{m}$, invariant under $\mathfrak{k}$. There exists a $\mathfrak{k}$-invariant supplement $\mathfrak{n}$ to $\mathfrak{r}$ in $\mathfrak{m}$, and $\mathfrak{n} \neq 0$. Then $B|\mathfrak{n}$ is non-degenerate, and we are back to d).

f) There is no proper $\mathfrak{k}$-invariant subspace $\mathfrak{m}$ of $\mathfrak{p}$ which is isotropic for $B$.

Let $\mathfrak{m}$ be one. If $\dim \mathfrak{m} < \dim \mathfrak{p}/2$, then the orthogonal subspace $\mathfrak{n}$ to $\mathfrak{m}$ has dimension $> \dim \mathfrak{p}/2$ and is invariant under $\mathfrak{k}$, and the restriction of $B$ to $\mathfrak{n}$ is non-zero. We are back to e). There remains to consider the case where $\dim \mathfrak{p}$ is even, $\dim \mathfrak{m} = \dim \mathfrak{p}/2$, $B|\mathfrak{p}$ is hyperbolic, and $\mathfrak{m}$ is maximal isotropic. Moreover, the representation of $\mathfrak{k}$ in $\mathfrak{m}$ is irreducible; otherwise we would be back to the case just treated. There exists a supplement $\mathfrak{n}$ to $\mathfrak{m}$ in $\mathfrak{p}$ which is $\mathfrak{k}$-invariant. Then, again, it has to be maximal isotropic. By c),

$$(4) \qquad [\mathfrak{m}, \, \mathfrak{m}] = [\mathfrak{n}, \, \mathfrak{n}] = 0;$$

hence by b)

$$(5) \qquad \mathfrak{k} = [\mathfrak{m}, \, \mathfrak{n}].$$

It follows from 1) and 4) that $\mathfrak{m}$ and $\mathfrak{n}$ consist of nilpotent matrices. The normalizer of $\mathfrak{n}$ in $\mathfrak{m}$ is $\mathfrak{k}$-invariant, hence reduced to zero in view of (5) and the irreducibility of the representation of $\mathfrak{k}$ in $\mathfrak{m}$. Therefore $\mathfrak{n}$ is the nilpotent radical of $\mathfrak{k} \oplus \mathfrak{n}$, and the latter is the normalizer of $\mathfrak{n}$. By Theorem 2 in [Bu1, 8, §10], $\mathfrak{k} \oplus \mathfrak{n}$ is parabolic. But $\mathfrak{k}$ is assumed to be semisimple, whence a contradiction. This concludes the proof.

*Remarks.* 1) In this paper, 8.0 is only needed in case $E = k_s$, $\mathfrak{g} = \mathfrak{g}_s$ and $\mathfrak{k} = \mathfrak{h}_s$. If it holds after extension of the groundfield, it is clearly already true in the original situation. Since $k_s$ may be embedded into $\mathbf{C}$, this reduces us to the case where $E = \mathbf{C}$, $\mathfrak{g} = \mathbf{sl}_n \mathbf{C}$ and $\mathfrak{k}$ is the Lie algebra of $\mathbf{SO}_n(\mathbf{C})$. A reader who would have a direct argument in that last case could then avoid any recourse to 8.0.

2) As the fixed point set of an involution, $\mathfrak{k}$ is always reductive. The assumption $\mathfrak{k}$ semisimple has been used only in the last step of the proof. Some restriction is necessary, since otherwise $\mathbf{sl}_2$ already provides a counterexample, with $\sigma$ having the diagonal matrices as fixed point set. Over $\mathbf{C}$, other counterexamples are given by the complexifications of hermitian symmetric pairs. In fact, the proof leads to a complete description of the cases where $\mathfrak{k}$ is not proper maximal, namely,

$$\mathfrak{g} = \mathfrak{k} \oplus \mathfrak{n} \oplus \mathfrak{m}$$

where $\mathfrak{n}$, $\mathfrak{m}$ are commutative, are the nilpotent radicals of two parabolic subalgebras with maximal reductive subalgebra $\mathfrak{k}$, the representation of $\mathfrak{k}$ in $\mathfrak{n}$ and $\mathfrak{m}$ are contragredient of one another, irreducible, the split center of $\mathfrak{k}$ is one-dimensional and acts by dilations on $\mathfrak{m}$ and $\mathfrak{n}$. Conversely, given such a decomposition of $\mathfrak{g}$, the map of $\mathfrak{g}$ onto itself which is the identity on $\mathfrak{k}$ and $-\mathrm{Id}$. on $\mathfrak{m} \oplus \mathfrak{n}$ is obviously an automorphism.

School of Mathematics, Institute for Advanced Study, Princeton, New Jersey
08540
   *E-mail address*: `borel@math.ias.edu`